\documentclass{article}
\usepackage{amssymb,amsfonts,amsmath,amsthm,subfigure}
\usepackage{authblk}
\usepackage{graphicx}
\numberwithin{equation}{section} 

\newtheorem{thm}{Theorem}[section]
\newtheorem{defn}[thm]{Definition}

\newtheorem{lemma}[thm]{Lemma}

\newtheorem{rmrk}[thm]{Remark}

\newcommand{\e}{\varepsilon}
\newcommand{\R}{\mathbb{R}}
\newtheorem{prop}[thm]{Proposition}

\newcommand{\abs}[1]{\left\vert{#1}\right\vert}

\newcommand{\C}{\mathbb{C}}
\newcommand{\ba}{\begin{array}}
\newcommand{\ea}{\end{array}}

\newcommand{\bthm}{\begin{thm}}
\newcommand{\ethm}{\end{thm}}
\newcommand{\bstp}{\begin{stp}}
\newcommand{\estp}{\end{stp}}
\newcommand{\blemma}{\begin{lemma}}
\newcommand{\elemma}{\end{lemma}}
\newcommand{\bprop}{\begin{prop}}
\newcommand{\eprop}{\end{prop}}
\newcommand{\bpf}{\begin{pf}}
\newcommand{\epf}{\end{pf}}
\newcommand{\bdefn}{\begin{defn}}
\newcommand{\edefn}{\end{defn}}
\newcommand{\brk}{\begin{rmrk}}
\newcommand{\erk}{\end{rmrk}}
\newcommand{\bcrl}{\begin{crl}}
\newcommand{\ecrl}{\end{crl}}
\newcommand{\norm}[1]{\left\|#1\right\|}

\newcommand{\beqn}{\begin{equation}}
\newcommand{\eeqn}{\end{equation}}

\renewcommand{\leq}{\leqslant}
\renewcommand{\geq}{\geqslant}
\newcommand{\da}{\dagger}
\newcommand{\rpl}{{\rm{Re}}\,\lambda_1}
\newcommand{\ipl}{{\rm{Im}}\,\lambda_1}
\newcommand{\mL}{\mathcal{L}}
\newcommand{\rec}{\mathcal{R}}
\newcommand{\intR}{\int_{\mathcal{R}}}
\newcommand{\beq}{\begin{equation}}
\newcommand{\eeq}{\end{equation}}
\newcommand{\bea}{\begin{eqnarray}}
\newcommand{\eea}{\end{eqnarray}}
\begin{document}
\renewcommand\Authfont{\small}
\renewcommand\Affilfont{\itshape\footnotesize}
\title{Kinematic and Dynamic Vortices in a Thin Film Driven by an Applied Current and Magnetic Field}

\author[1]{Lydia Peres Hari\footnote{ lydia@fermat.technion.ac.il}}
\author[1]{Jacob Rubinstein\footnote{koby@techunix.technion.ac.il}}
\author[2]{Peter Sternberg\footnote{Corresponding author, sternber@indiana.edu}}
\affil[1]{Department of Mathematics, Israel Institute of Technology, Haifa 32000, Israel}
\affil[2]{Department of Mathematics, Indiana University, Bloomington, IN 47405}


\maketitle
\begin{abstract}
Using a Ginzburg-Landau model, we study the vortex behavior of a rectangular thin film superconductor subjected to an applied
current fed into a portion of the sides and an applied magnetic field directed orthogonal to the film.
 Through a center manifold reduction we develop a rigorous bifurcation theory for the appearance of periodic solutions
in certain parameter regimes near the normal state. The leading order dynamics yield in particular a motion law for
kinematic vortices moving up and down the center line of the sample. We also present computations that reveal the co-existence
and periodic evolution of kinematic and magnetic vortices.
\end{abstract}

\section{Introduction}

We consider a thin superconducting sample occupying a rectangle. The sample is subjected to normal electric current that enters through a lead on one side and leaves through another lead at the opposite side. In addition the sample is subjected to a magnetic field oriented in the direction perpendicular to the rectangle's plane. Our main interest is in the regime in the parameter space where different physical quantities, such as the total current in the sample and the order parameter, are time-periodic. We pay special attention to the formation and motion of vortices in the sample.

The problem above is a natural generalization of the simpler case of a finite superconducting one-dimensional wire subjected to normal current that is fed into one of its endpoints. This problem received considerable attention since it is a canonical case of co-existence of normal current and superconducting current. Moreover, it is known that in this setting there exists a regime of prescribed current $I$ and temperature $T$, where the sample's behavior is time-periodic. In addition, in this regime the superconducting order parameter $\psi(x,t)$ vanishes at the wire's center at specific points in time that are separated by a fixed period. Such zeros of $\psi$ are called phase slip centers (PSC). These phenomena and others are studied, numerically and experimentally, by a number of authors including \cite{krba,krwa,laam}, and many of the results are summarized by Ivlev and Kopnin \cite{ivko}.

A recent study in \cite{rsm} and \cite{rsz} presents a comprehensive theory that explains the different patterns observed in this wire setting. The key idea is that the underlying system of equations enjoys a PT symmetry, that is, symmetry under
complex conjugation and the transformation $x \rightarrow -x$. This symmetry enabled the authors to perform a rigorous bifurcation study of the problem, and to deduce that, under certain conditions, the order parameter bifurcates, as the temperature is lowered beyond a critical value, from the normal state $\psi \equiv 0$ to a nontrivial state. Moreover, this bifurcation is shown to be of Hopf type, and this explains the periodic nature of the solution and the periodic appearance of isolated zeros of $\psi$.

In the present study we look at a more realistic geometry of a finite strip. Moreover,  we consider not just forced electric current, but also the effect of an external magnetic field. The problem is analyzed numerically in \cite{bmp}. The authors observe, just as in the one-dimensional setting, a periodic behavior of a number of physical quantities, including periodic appearance and motion of vortices. Therefore, our goal is to derive a theory that explains the observed patterns and vortex motion.

The issue of vortices is of particular interest.  They are defined as isolated zeros of the order parameter, and they are characterized by their topological degree in the $(x,y)$ plane. The appearance of vortices in superconducting samples subjected to an applied magnetic field is of course well-known. Therefore, we expect them in our setting even in the absence of forced electric currents. What makes the present problem interesting is that also the opposite is true; namely, vortices form, for appropriate range of values of $I$ and $T$, even if no magnetic field is applied. Therefore, one can classify the vortices here into magnetic vortices generated by the applied magnetic field in the absence of any applied current, and kinematic vortices generated, as will be shown below, by the forced electric field and the special symmetry of the problem in the absence of any magnetic field.

When both applied magnetic and electric fields are present, it is not as
clear how to make a distinction between the two kinds of vortices. As will be shown below,
the symmetry of the problem implies that some vortices are formed and move time-periodically on the center line of the rectangle for large enough $I$ and for a range of applied magnetic fields $h$. We term them kinematic vortices. In certain cases, as presented below, these kinematic vortices collide and move off the center line line. We term such vortices,
``born" from kinematic vortices, `kinematic' as well.

Other vortex phenomena, not usually observed in more standard Ginzburg-Landau settings, include the time-periodic emergence of vortex pairs that are of opposite degrees, which we term vortex/anti-vortex pairs: typically stable vortex configurations of Ginzburg-Landau vortices only involve vortices of the same degree but in this periodic setting that is not always the case.  We also show that for some range of $I$ and $h$ vortices of {\em the same degree} move towards each other and collide before moving away from each other. Again, this is a process that is atypical to more familiar Ginzburg-Landau settings.

The analysis of the present problem follows to some extent the lines of the one-dimensional wire problem. Namely, we construct a proper framework that enables us to use the Center Manifold Theorem to study the bifurcation picture, and thereby establish the existence of a Hopf bifurcation. Moreover, in both cases a key factor is played by the spectrum of the underlying linear Schr\"odinger operator. However, there are a few important differences between the one-dimensional case and the problem considered here. First, the construction of the center manifold requires certain apriori estimates on the solutions to the underlying differential equations. These estimates are harder to obtain in the two-dimensional setting. A second important difference relates to the vortex motion. While in the one-dimensional case the PSCs just appear momentarily at a fixed periodic sequence of instances, the kinematic vortices in the two-dimensional problem are present for periodic finite time inte
 rvals. Moreover, they tend to move along or near the $y$-axis (which is the center line of the rectangle) as we shall show.

It is interesting to note that implications of PT-symmetry seem to appear in a number of quite different physical problems. For instance, we mention applications to quantum mechanics \cite{bebo}, \cite{cdv}, to hydrodynamic instability \cite{shk1}, \cite{shk2}, and to optics \cite{zia1}, \cite{zia2}. We also mention several recent rigorous studies within Ginzburg-Landau theory
that incorporate magnetic effects along with applied currents in a variety of asymptotic regimes, \cite{Almog,AHP,dwz,ST,Tice}. One aspect of our investigation that distinguishes it from others, however, is that it is to our knowledge the first to capture a motion law for Ginzburg-Landau vortices that is not based on the assumption of large Ginzburg-Landau parameter.

In the next section we formulate the problem and the underlying equations. The bifurcation analysis is performed in section 3. In particular we establish there the existence of a center manifold for a certain regime in the $(I,T)$ plane. In section 4 we consider the formation of vortices and their motion and discuss some computational work on the problem.

\label{intro}
\section{Formulation of Problem}
\begin{figure}
\begin{tabular}{c}
\begin{tabular}{cc}
\includegraphics[width=0.9\textwidth]{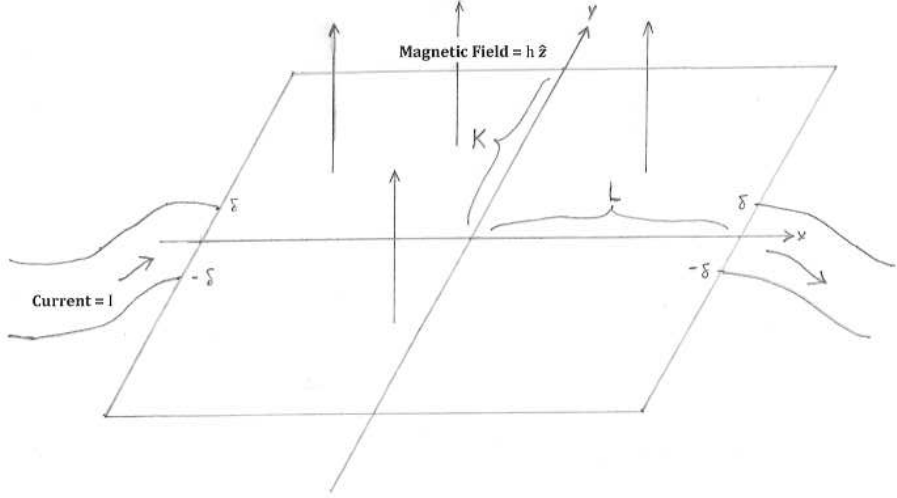}
\end{tabular}
\end{tabular}
\caption{A thin film superconductor subjected to applied current and magnetic field.}
\label{rectangle}
\end{figure}

We consider a superconducting material occupying a thin rectangular
box with dimensions $-L<x<L,\;-K<y<K$ and say $0<z<\eta$ where
$\eta$ is assumed to be much smaller than the coherence length or
penetration depth, allowing us to work within the thin film 2d
approximation of Ginzburg-Landau. In this approximation, we take the
complex-valued order parameter $\Psi=\Psi(x,y,t)$ and the
real-valued electric potential $\phi=\phi(x,y,t)$ to be defined on
$\rec\times[0,\infty)$ where we denote $\rec:=[-L,L]\times[-K,K]$
and we ignore any induced magnetic field. Within this rectangular
geometry, we assume the presence of leads forcing in electric
current of magnitude $I$ through the sides $x=\pm L$, along the
subinterval $-\delta<y<\delta$ for some positive $\delta<K.$
Additionally we assume that the thin film is subjected to an applied
magnetic field of size $h$ oriented perpendicular to the rectangular
cross-section. See Figure \ref{rectangle}.

Introducing the applied magnetic potential
$A_0:=(-y,0)$, this applied field is then given by $h\nabla\times
A_0$ and the thin film limit of Ginzburg-Landau takes the nondimensionalized form
\begin{eqnarray}
& \Psi_t+i\phi\Psi=\left(\nabla-ihA_0\right)^2\Psi+(\Gamma-\abs{\Psi}^2)\Psi\quad\mbox{for}\;(x,y)\in\rec,\;t>0,\label{psieqn}\\
&\Delta \phi=\nabla\cdot\bigg(\frac{i}{2}\{\Psi\nabla\Psi^*-\Psi^*\nabla \Psi\}-\abs{\Psi}^2hA_0\bigg)\quad\mbox{for}\;(x,y)\in\rec,\;t>0,\label{phieqn}
\end{eqnarray}
subject to the boundary conditions
\begin{eqnarray}
&\Psi(\pm L,y,t)=0\;\mbox{for}\;\abs{y}<\delta,\label{bc1}\\
&\Psi_x(\pm L,y,t)+ihy\Psi(\pm L,y,t)=0\;\mbox{for}\;\delta<\abs{y}\leq K,\label{bc2}\\
&\Psi_y(x,\pm K,t)=0\;\mbox{for}\;\abs{x}\leq L,\label{bc4}\\
&\phi_x(\pm L,y,t)=\left\{\begin{matrix}-I&\;\mbox{for}\;\abs{y}<\delta,\\ 0&\;\mbox{for}\;\delta<\abs{y}<K\end{matrix}\right.\label{bc5}\\
&\phi_y(x,\pm K,t)=0\;\mbox{for}\;\abs{x}\leq L.\label{bc6}
\end{eqnarray}
along with the initial condition $\Psi(x,y,0)=\Psi_{init}(x,y).$ As a
convenient normalization we take
$\intR\phi=0.$ The system of equations above are collectively known as the Time-Dependent Ginzburg-Landau (TDGL) model.

The parameter
$\Gamma$ is proportional to $T-T_c$ where $T$ is temperature and
$T_c$ is the critical temperature below which the normal (zero)
state loses stability in the absence of any applied fields. Note
that \eqref{phieqn} is simply the requirement of conservation of
total (normal $+$ superconducting) current. Also, we remark that the
boundary conditions \eqref{bc2}-\eqref{bc4} are the standard
superconductor/vacuum conditions on $\Psi$, while $\eqref{bc1}$
reflects the presence of the normal leads.

There are numerous investigations of current-driven superconducting
wires and thin films that utilize a model based on Ginzburg-Landau theory,
and these have generally reported reasonable agreement between
theory and experiment. Regarding \eqref{psieqn}, we should mention
that some studies on this problem such as \cite{bepv,bmp,mmpvp} replace
the left-hand side with the modification
\[\frac{u}{\sqrt{1+\gamma^2\abs{\Psi}^2}}\bigg(\frac{\partial}{\partial
t}+i\phi+\frac{\gamma^2}{2}\frac{\partial\abs{\Psi}^2}{\partial
t}\bigg)\Psi,\] where $u$ and $\gamma$ are material parameters.
However, we believe that the standard and simpler evolution equation
\eqref{psieqn} corresponding to the choices $u=1$ and $\gamma=0$
captures the main features of the problem, an opinion shared by the
authors of \cite{agkns} whose computational comparisons suggest that the
modification does not have a large effect in this setting. What is
more, in what follows, we shall concentrate on bifurcation from
the normal state $\Psi\equiv 0$, so the smallness of the amplitude
should mute the effect of this modification even more.

\section{Analysis of the model}

For a given $\Psi$, and a given value $h\geq 0$, let us decompose
the solution $\phi$ to \eqref{phieqn}, \eqref{bc5}, \eqref{bc6} as
\[\phi=I\phi^0+\tilde{\phi}\]  and where $\phi^0$ is harmonic and satisfies
the boundary conditions
\begin{eqnarray}
&\phi^0_x(\pm L,y)=\left\{\begin{matrix}-1&\;\mbox{for}\;\abs{y}<\delta,\\ 0&\;\mbox{for}\;\delta<\abs{y}<K,\end{matrix}\right.\label{bcx}\\
&\phi^0_y(x,\pm K)=0\;\mbox{for}\;\abs{x}\leq L,\label{bcy}
\end{eqnarray} while $\tilde{\phi}$ satisfies
\eqref{phieqn} subject to homogeneous Neumann boundary conditions on
all portions of the rectangular boundary. We note that
$\tilde{\phi}$ but depends
on $\Psi$ so we will often write $\tilde{\phi}$ as
$\tilde{\phi}[\Psi]$ to emphasize this dependence.

We will occasionally make use of the properties
\begin{equation}
\phi^0(-x,y)=-\phi^0(x,y)\quad\mbox{and}\quad \phi^0(x,-y)=\phi^0(x,y),\label{oddeven}
\end{equation}
which are easy to check.

The normal state in this setting corresponds to $\Psi\equiv 0$ and
$\phi\equiv I\phi^0$. We will pursue a bifurcation analysis about
this normal state, and therefore a crucial role will be played by the
linear eigenvalue problem
\begin{equation}
 \mL[u]:=\big(\nabla -ih A_0)^2u-iI\phi^0 u=-\lambda u\quad\mbox{for}\;\abs{x}<L,\;\abs{y}<K,\label{evpm}
\end{equation}
subject to the boundary conditions \eqref{bc1}--\eqref{bc4}. Though
we do not indicate it in our notation, it is understood that
$\mL$ and therefore all of its eigenvalues and eigenfunctions
depend on the parameters $L,K,\delta,h$ and $I$. We summarize below
the key properties of the corresponding eigenvalues and
eigenfunctions that will be needed in the analysis to follow.
\blemma\label{evproblem} The spectrum of $\mL$ consists only of
point spectrum, denoted by $\{\lambda_j\}$ with corresponding
eigenfunctions $\{u_j\}.$ If $(\lambda_j,u_j)$ is
an eigenpair satisfying \eqref{evpm} then
\begin{equation} {\rm{Re}}\,\lambda_j>0,\quad\mbox{and}\quad
\abs{{\rm{Im}}\,\lambda_j}<
\norm{\phi^0}_{L^{\infty}}I.\label{ReIm}\end{equation} Thus, in
particular we may order the eigenvalues
$\lambda_1,\lambda_2,\ldots$ according to the size of
their real part, with $0<{\rm{Re}}\,\lambda_1\leq
{\rm{Re}}\,\lambda_2\leq \ldots.$ The PT-symmetry of the
operator is reflected in the fact that if
$(\lambda_j,u_j)$ is an eigenpair then so is
$(\lambda^*_j,u_j^{\da})$ where
$u_j^{\da}(x,y):=u_j^*(-x,y).$ \elemma

When working with $u_1$and $u_2$, we will choose the normalization
\begin{equation}
\intR u_1^2=1=\intR u_2^2. \label{normal1}
\end{equation}
Also, as a matter of convention, when $\lambda_1$ is non-real, we associate $u_1$ with the eigenvalue $\lambda_1$ having positive imaginary part,
and $u_2$ with $\lambda_2^*$.

\begin{proof} The fact that the spectrum consists solely of
eigenvalues follows from standard compact operator theory. The
conditions in \eqref{ReIm} follow from multiplication of the
equation $\mL[u]=-\lambda u$ by $u^*$ and integration over the
rectangle. This leads to the identity
\begin{equation}
\lambda=\frac{\intR\abs{\left(\nabla \pm i h A_0\right) u}^2}{\intR\abs{u}^2}+
iI\frac{\intR \phi^0 \abs{u}^2}{\intR \abs{u}^2},\label{lamid}
\end{equation}
implying \eqref{ReIm}.

The final claim of the lemma follows by noting that
$\mL[u^{\da}]=\mL[u].$
\end{proof}

\begin{figure}
\subfigure[]{\label{aaa1}
\includegraphics[width=6cm]{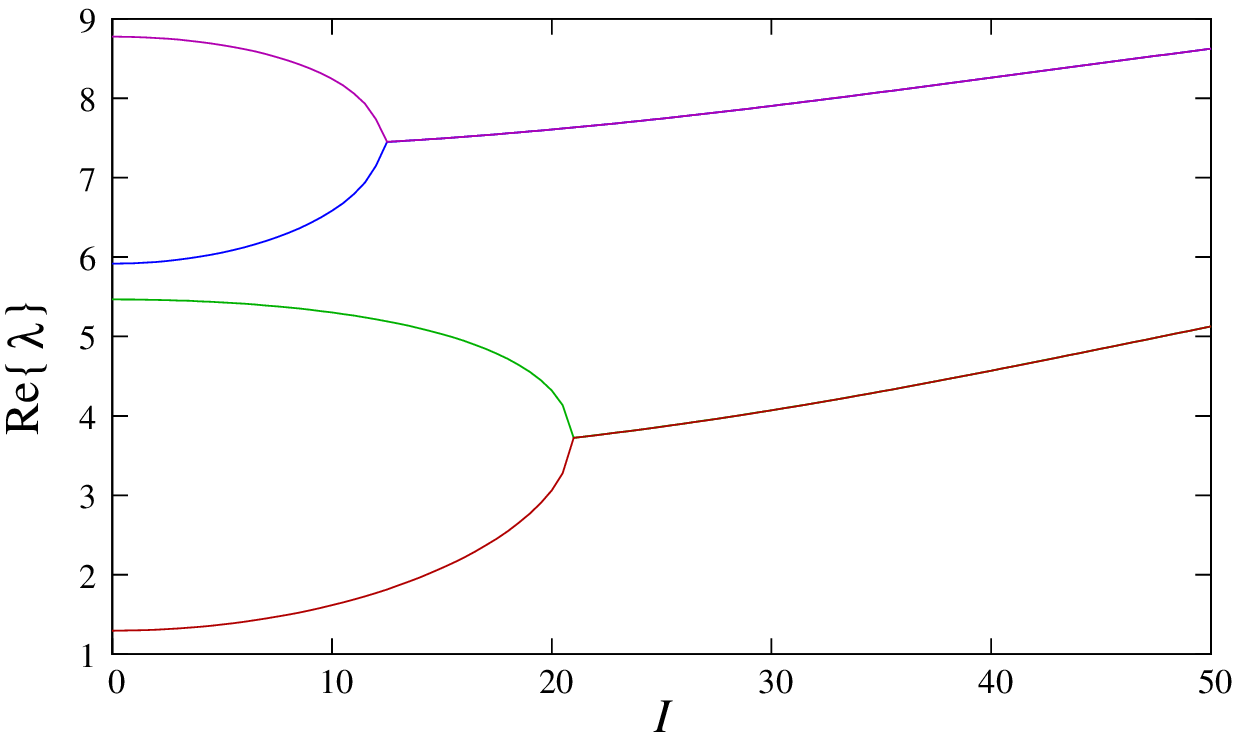}
}
\subfigure[]{\label{aaa2}
\includegraphics[width=6cm]{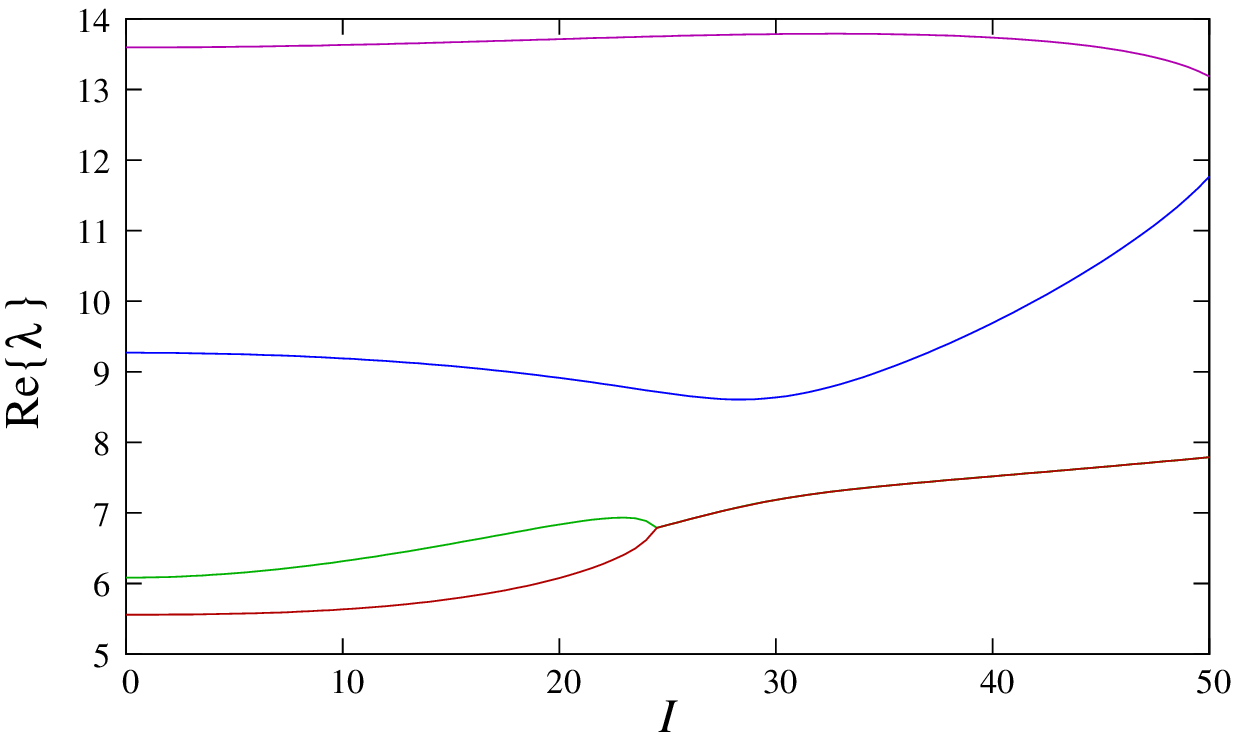}
}
\caption{The real part of the spectrum of $\mL$. The first 4 eigenvalues
are drawn for the parameters $L=1, K=2/3, \delta = 1/6$. In (a) $h=0$, while in (b) $h=7.5$} \label{eig_val}
\end{figure}

A numerical analysis of the eigenvalue problem \eqref{evpm}
indicates that, fixing all other parameters (i.e. $K,L,\delta$ and
$h$) and varying only the current $I$, there exists a critical value
$I_c$ depending on these other parameters, such that
$\lambda_1$ is real for $I\leq I_c$ and $\lambda_1$ is
non-real for $I>I_c$. This complexification of $\lambda_1$ arises through
a collision with another eigenvalue and the two eigenvalues emerge from this
collision, that is for $I>I_c$, as complex conjugate pairs.
Examples of eigenvalue collisions are shown in Figures \ref{eig_val}a and \ref{eig_val}b.
In both cases the simulated geometry is $L=1, K=2/3, \delta = 4/15$. In Figure \ref{eig_val}a we took $h=0$, while in Figure \ref{eig_val}b
we took $h=7.5$.

We note in passing that in light of
\eqref{oddeven} and \eqref{lamid}, the evenness of $\abs{u_j}$
evidently implies reality of $\lambda_j$ so the passage to
complex eigenvalues with increased values of $I$ carries with it a
certain symmetry breaking of the corresponding eigenfunction. This
phenomenon has been established rigorously for the one-dimensional
version of \eqref{evpm}--that is, for the problem depending only on
$x$ and with $h=0$ in \cite{shk1,shk2}. The reality of the spectrum
for $I$ positive and sufficiently small in this two-dimensional setting should follow by the type of
perturbation analysis to be found in \cite{cgs}, since for $I=0$,
the spectrum is clearly real, cf. \eqref{ReIm}. However, since here
we are interested in the regime where the first eigenvalue is
complex, we do not pursue this point further.

\begin{figure}
\begin{tabular}{c}
\begin{tabular}{cc}
\includegraphics[width=0.5\textwidth]{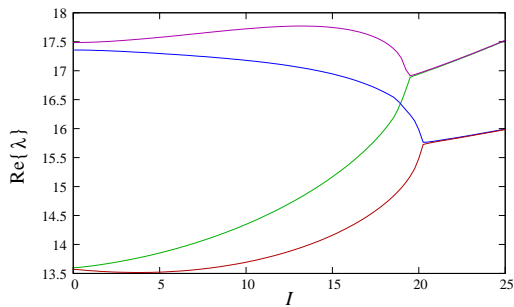}
\end{tabular}
\end{tabular}
\caption{The real parts of the four leading eigenvalues for the parameter values $L=1, K=2/3, \delta = 4/15, h=20, I=25$. Note how $\lambda_3$
passes through $\lambda_2$ and then collides with $\lambda_1$ at $I=I_c\approx 20$. } \label{var_beta1}
\end{figure}

Regarding eigenvalue collisions, we also wish to note another new phenomenon
for this two-dimensional problem with the incorporation of magnetic effects that is
not observed for the one-dimensional problem. Our computations reveal that for certain
large enough values of applied magnetic field $h$, as $I$ increases through the regime $0<I<I_c$,
the second and third (still real) eigenvalues pass through each other and it is ultimately
what was originally labeled as the third eigenvalue that collides with $\lambda_1$ at $I_c$.
We say the third eigenvalue ``passes through" the second rather than ``collides" with it because
as $I$ varies through the point where $\lambda_2=\lambda_3$, both eigenvalues remain simple
and the corresponding eigenfunctions vary smoothly without incident.
This for example
is the scenario for the parameter values
$L=1, K=2/3, \delta = 4/15, h=20, I=25$. See Figure \ref{var_beta1}. We will return to this set of
parameter values at the end of the article to discuss anomalous vortex behavior as well.

At this point we fix any values of $K,L,\delta,h$ and then pick $I$
sufficiently large so that ${\rm{Im}}\,\lambda_1\not= 0$. We then
consider our problem \eqref{psieqn}--\eqref{bc6} with $\Gamma$ given
by $\Gamma={\rm{Re}}\,\lambda_1+\e$ where $\e$ is small and
positive. Note that given the temperature dependence of $\Gamma$,
this corresponds to lowering the temperature just below the value
where the normal states becomes linearly unstable.

We now introduce the linear operator $\mL_1[u]:=\mL[u]+(\rpl)\, u$.
Recalling Lemma \ref{evproblem} in the scenario where $\lambda_1$ is not real, the
corresponding first eigenfunctions, say $u_1$ and $u_2$, satisfy
\begin{equation}
u_2=u_1^{\da}\quad \mbox{with}\quad\mL_1[u_1]=-i\,\ipl\,u_1,\quad \mL_1[u_2]=i\,\ipl\,u_2.
\label{u1u2}\end{equation}
For later use, we also note that if we introduce the adjoint operator $\mL_1^{\star}$ satisfying
\[
\intR \mL_1[u]v=\intR \mL_1^{\star}[v]u,\]
then we can readily identify $\mL_1^{\star}[v]$ as simply the operator
given by $\big(\nabla +ih A_0)^2v-iI\phi^0v+(\rpl)\,v$. Furthermore, the
fact that $\phi_0$ is even in $y$ (cf. \eqref{oddeven}) reveals that
the functions $u^{\star}_j(x,y):=u_j(x,-y)$ for $j=1,2$ satisfy the equations
\[
\mL_1^{\star}[u^{\star}_1]=-i\,\ipl\,u^{\star}_1,\quad  \mL_1^{\star}[u^{\star}_2]=i\,\ipl\,u^{\star}_2.\]
Consequently, we see that
\begin{equation*}
-i\,\ipl\intR u^{\star}_1u_2=\intR \mL_1^{\star}[u^{\star}_1]u_2=\intR \mL_1[u_2]u^{\star}_1=i\ipl\intR u^{\star}_1u_2,\end{equation*}
with a similar relation holding between $u^{\star}_2$ and $u_1$. Hence,
$\intR  u^{\star}_1u_2=0=\intR u^{\star}_2u_1.$
Clearly, this orthogonality holds between any two eigenfunctions $u_j$ and $u_k$ of the operator
$\mL_1$ corresponding to distinct eigenvalues, namely
\begin{equation}
\intR  u^{\star}_ju_k=0\quad\mbox{for}\;j\not=k.\label{orthog}
\end{equation}

Of course, in the special case of no applied magnetic field, i.e. $h=0$, we have
$\mL=\mL^{\star}$ and $u_1=u^{\star}_1$ and through \eqref{oddeven} we see then that each eigenfunction is
even in $y$.

Let us now return to
the \eqref{psieqn}-\eqref{phieqn} with the choice  $\Gamma={\rm{Re}}\,\lambda_1+\e$ and re-express the system as
a single nonlinear, non-local equation for the complex-valued order parameter $\Psi$:
\begin{equation}
\Psi_t=\mL_1[\Psi]+\e \Psi+\mathcal{N}(\Psi),\label{pde}\end{equation} where
\begin{equation}
\mathcal{N}(\Psi):=-\abs{\Psi}^2\Psi-i\tilde{\phi}[\Psi]\Psi.\label{nonlinear}
\end{equation}

When expressed in this form, the problem can be rigorously solved
for $\e\ll 1$ via a center manifold reduction. The analysis is
similar to that carried out for the one-dimensional (thin wire)
problem in \cite{rsz}, Prop. 6.8, so we will only mention the key
steps and those parts of the calculation where there are changes.

As regards the linear part of \eqref{pde}, the key point is that the
operator $-\mL_1$ is sectorial, in light of \eqref{ReIm}, cf. \cite{Henry}.
Regarding estimates on the cubic, nonlocal nonlinearity
$\mathcal{N}$, the analysis differs from that in \cite{rsz} in
that for the one-dimensional problem it is easy to check that
$\mathcal{N}$ is a bounded map from $H^1$ to $H^1$ (cf. \cite{rsz},
Lemma 6.3), while in the present two-dimensional setting this is no
longer true--it just barely misses. One way to overcome this obstacle is
by viewing $\mathcal{N}$ as a mapping from an
interpolation space between $H^1$ and $H^2$ into $L^2$. The details
of pursuing this strategy can be found in section 3.4 of \cite{mw},
where the authors execute a center manifold construction relevant to
bifurcation from the normal state without applied electric current.
Alternatively, one can view $\mathcal{N}$ as a map from $H^2(\mathcal{R})$ into
$H^1(\mathcal{R})$. We describe how to make the necessary estimates for this latter
approach.

Writing $\mathcal{N}=\mathcal{N}_1+\mathcal{N}_2$ with $\mathcal{N}_1(\Psi):=-\abs{\Psi}^2\Psi$ and
$\mathcal{N}_2(\Psi):=-i\tilde{\phi}[\Psi]\Psi,$ it is an easy application of H\"older's inequality
to make the estimate
\beq
\norm{\mathcal{N}_1(\Psi)}_{H^1(\mathcal{R})}\leq C\norm{\Psi}^3_{H^2(\mathcal{R})}.
\label{N1est}
\eeq
To make a similar estimate on $\mathcal{N}_2$, let us first write the PDE coming from \eqref{phieqn} for
$\tilde{\phi}(\Psi)$ as
\beq
\Delta \tilde{\phi}=\nabla\cdot \bigg(j(\Psi)-\abs{\Psi}^2A_0\bigg),\label{newtilde}\eeq
subject to homogeneous Neumann conditions and the normalization condition of mean zero,
where we have introduced the notation $j(\Psi):=\frac{i}{2}\{\Psi\nabla\Psi^*-\Psi^*\nabla \Psi\}.$
Fixing any $p\in (1,2)$ we may use H\"older's inequality and Sobolev imbedding to make the estimate
\begin{eqnarray*}
&&
\intR \abs{Dj(\Psi)}^p\leq C\intR \bigg(\abs{\Psi}^p\abs{D^2\Psi}^p+\abs{D\Psi}^{2p}\bigg)\\
 && \leq C\bigg(\norm{\Psi}^p_{H^2(\mathcal{R})}\norm{\Psi}^p_{L^{2p/(2-p)}(\mathcal{R})}+\norm{\Psi}^{2p}_{W^{1,2p}(\mathcal{R})}\bigg)\leq
 C\norm{\Psi}^{2p}_{H^2(\mathcal{R})}.
\end{eqnarray*}
Hence, $\norm{Dj(\Psi)}_{L^p(\mathcal{R})}\leq C\norm{\Psi}^2_{H^2(\mathcal{R})}.$ It is also easy to estimate
\[\norm{j(\Psi)}_{L^p(\mathcal{R})}+\norm{\abs{\Psi}^2A_0}_{W^{1,p}(\mathcal{R})}\leq
C\norm{\Psi}^2_{H^2(\mathcal{R})},\] so
we conclude that \[\norm{\nabla\cdot\big(j(\Psi)-\abs{\Psi}^2A_0\big)}_{L^{p}(\mathcal{R})}\leq C\norm{\Psi}^2_{H^2(\mathcal{R})}.\]
Then appealing to the Calderon-Zygmund inequality (cf. \cite{W}, Chapter 2), equation \eqref{newtilde} implies that
\[
\norm{\tilde{\phi}(\Psi)}_{W^{2,p}(\mathcal{R})}\leq C\norm{\Psi}^2_{H^2(\mathcal{R})}.\]
From the Sobolev imbedding theorem and H\"older's inequality it then follows easily that
$\norm{\mathcal{N}_2(\Psi)}_{H^1(\mathcal{R})}\leq C\norm{\Psi}^3_{H^2(\mathcal{R})}$.
Combining this last estimate with \eqref{N1est} we arrive at the estimate on the nonlinearity:
\beq
\norm{\mathcal{N}(\Psi)}_{H^1(\mathcal{R})}\leq C\norm{\Psi}^3_{H^2(\mathcal{R})}.\label{nonlinearest}
\eeq
The upshot is that for each small $\e$, one can construct a center
manifold $\mathcal{M}_\e$ as a graph $v\mapsto\Phi(v,\e)$ in $H^2({\mathcal{R}})$ over
the center subspace $\mathcal{S}:={\rm{span}}\,\{u_1,u_2\}$, applying for example, the version
of the Center Manifold Theorem to be found in \cite{HI}, Theorem 2.9. More precisely, there exist positive constants
$\delta_0$ and $\e_0$, such that for any $\e$ satisfying
$\abs{\e}\leq \e_0$ one can define $
\mathcal{M}_{\e}:=\{\Phi(v,\e):\,v\in
\mathcal{S},\;\norm{v}_{H^2(\mathcal{R})}<\delta_0\},$ and $\mathcal{M}_\e$
enjoys the following properties:\\

\noindent (i) The center manifold is locally invariant for the flow
\eqref{pde} in the sense that if $\abs{\e}<\e_0$ and the initial
data $\psi_0$ lies on $\mathcal{M}_\e$, then so does the solution
$\psi^{\e}$ to \eqref{pde} so long as
$\norm{\psi^\e(\cdot,t)}_{H^2(\mathcal{R})}$ stays sufficiently small.
Hence, for such initial data, one can describe the resulting
solution $\psi^{\e}(t)=\psi^{\e}(\cdot,t)$ through two maps
$\alpha^{\e}_1,\,\alpha^{\e}_2:[0,\infty)\to\C$ via
$\psi^{\e}(t)=\Phi(\alpha^{\e}_1(t)u_1+\alpha^{\e}_2(t)u_2,\e).$ Since
$u_1$ and $u_2$ are fixed throughout, we will often write simply
\begin{equation}\psi^{\e}(t)=\Phi(\alpha^{\e}_1(t),\alpha^{\e}_2(t),\e)\label{cmr}\end{equation} for
the sake of brevity.

\noindent (ii) The center manifold is PT-symmetric, i.e.
$\psi\in\mathcal{M}_\e\implies\;\psi^{\da}\in \mathcal{M}_\e$.
Furthermore, if $v=v^{\da}$, then $\Phi(v,\e)=\Phi^{\da}(v,\e).$
Since $u_2=u_1^{\da}$, this implies
that if $\alpha^{\e}_2(0)=\alpha^{\e}_1(0)^*$,
and one solves \eqref{pde} subject to the PT-symmetric initial
conditions $\psi^{\e}(0)=\Phi(\alpha^{\e}_1(0),\alpha^{\e}_2(0),\e)$,
then the resulting functions
$\alpha^{\e}_1(t)$ and $\alpha^{\e}_2(t)$ describing the solution $\psi^{\e}$ at
any positive time $t$ will remain complex conjugates. \\

\noindent (iii) $\mathcal{M}$ contains all nearby bounded solutions
of \eqref{pde}, and in particular, it contains any nearby
steady-state or time-periodic solutions.\\

\noindent (iv)  Through an appeal to \eqref{nonlinearest}, the discrepancy between the center manifold and the
center subspace can be expressed through the estimate \beq
\norm{\Phi(v,\e)-v}_{H^2(\mathcal{R})} \leq
C_1\left(\norm{v}_{H^2(\mathcal{R})}^3+\abs{\e}\norm{v}_{H^2(\mathcal{R})}\right)\label{verytang}
\eeq which holds for any pair $(v,\e)$ such that $v\in \mathcal{S}$
with $\norm{v}_{H^2(\mathcal{R})}<\delta_0$ and $\abs{\e}<\e_0$, where $C_1$
is a positive constant independent of $v$ and $\e$.

Armed with these properties of the center manifold, one can then fix
any sufficiently small complex numbers $\alpha^{\e}_1(0)$ and
$\alpha^{\e}_2(0)$ as in (i) above, solve \eqref{pde} and then use the projection $\Pi_c$
onto the center subspace $\mathcal{S}$ to
obtain a reduced system of O.D.E.'s governing the evolution of $\alpha^{\e}_1$ and $\alpha^{\e}_2.$
In light of \eqref{orthog} we can write an arbitrary function $f$ as $f=c_1u_1+c_2u_2+u^{\perp}$ where
$u^{\perp}\in \big({\rm{span}}\,\{u^{\star}_1,u^{\star}_2\}\big)^{\perp}$. Hence the projection onto the center subspace
of an arbitrary function $f$ is given by
\[
\Pi_c(f)=\bigg(\frac{\intR u^{\star}_1f}{\intR u^{\star}_1u_1}\bigg)u_1+\bigg(\frac{\intR u^{\star}_2f}{\intR u^{\star}_2u_2}\bigg)u_2.\]
Thus, substituting the reduction \eqref{cmr} into \eqref{pde} and projecting, we obtain the following system through
the use of \eqref{u1u2}:
\bea &&\dot{\alpha}^{\e}_1u_1+\dot{\alpha}^{\e}_2u_2=(\e-i\,\ipl)\alpha^{\e}_1u_1+(\e+i\,\ipl)\alpha^{\e}_2u_2\nonumber\\
&&+\Pi_c\bigg(\mathcal{N}(\alpha^{\e}_1u_1+\alpha^{\e}_2 u_2)\bigg)+\bigg\{\Pi_c\bigg(\mathcal{N}(\Phi(\alpha^{\e}_1,\alpha^{\e}_2,\e))\bigg)
-\Pi_c\bigg(\mathcal{N}(\alpha^{\e}_1u_1+\alpha^{\e}_2 u_2)\bigg)\bigg\},\nonumber\\
\label{reduced}
\eea
where $\dot{}$ denotes a time derivative.

 Invoking \eqref{verytang}, one
 finds that the last expression above involving the difference of nonlinear terms is lower order so one is justified
 in initially ignoring it, solving the resulting simplified
 system of O.D.E.'s and then arguing that the behavior of solutions
 persists for the full system \eqref{reduced}. Again the details of this portion of the
 argument can be found in \cite{rsz}. In an abuse of notation, we will persist in using the
 notation $\alpha^{\e}_j$ to denote the solution to the truncated system in which the last expression is dropped.

Now using \eqref{orthog}, we integrate first against $u^{\star}_1$ and then against $u^{\star}_2$ to arrive at the system:
\bea
&&\dot{\alpha}^{\e}_1=(\e-i\,\ipl)\alpha^{\e}_1
+\frac{\intR u^{\star}_1\,\mathcal{N}(\alpha^{\e}_1u_1+\alpha^{\e}_2 u_2)}{\intR u^{\star}_1u_1},\label{redone}\\
&&\dot{\alpha}^{\e}_2=(\e+i\,\ipl)\alpha^{\e}_2
+\frac{\intR u^{\star}_2\,\mathcal{N}(\alpha^{\e}_1u_1+\alpha^{\e}_2 u_2)}{\intR u^{\star}_2u_2}.\label{redtwo}
\eea

In order to argue that \eqref{redone}-\eqref{redtwo} exhibits a Hopf bifurcation to a periodic state, we now restrict the flow to the PT symmetric
portion of the center subspace, and hence, in light of item (ii) above, to the PT symmetric subset of the center manifold. This amounts to the restriction
$\alpha^{\e}_2=(\alpha^{\e}_1)^*$ and allows us to only work with \eqref{redone}. At this juncture, we make a change of variables of the form
\[
a^{\e}:=\alpha^{\e}_1+c_1(\alpha^{\e}_1)^3+
c_2\abs{\alpha^{\e}_1}^2(\alpha^{\e}_1)^*+c_3((\alpha^{\e}_1)^*)^3,
\]
in order to convert the problem to its normal form, cf. \cite{HI}, chapter 3. Writing the cubic nonlinear term in \eqref{redone} as
\[
n_1(\alpha^{\e}_1)^3+n_2\abs{\alpha^{\e}_1}^2(\alpha^{\e}_1)^*+
n_3((\alpha^{\e}_1)^*)^3+n_4\abs{\alpha^{\e}_1}^2\alpha^{\e}_1,
\]
a tedious but direct calculation yields that with the choices $c_1=\frac{n_1}{2\ipl}\,i,\;c_2=-\frac{n_2}{2\ipl}\,i,$ and $c_3=-\frac{n_3}{4\ipl}\,i$,
the new variable $a^{\e}$ satisfies the simpler differential equation
\begin{equation}
\dot{a}^{\e}=(\e-i\,\ipl)a^{\e}+n_4\abs{a^{\e}}^2 a^{\e},\label{aeqn}
\end{equation}
with the coefficient $n_4$ given by
\beq
n_4=\frac{-
\intR\big(\abs{u_1}^2 u_1u^{\star}_1+2\abs{u_2}^2 u_1u^{\star}_1\big)-i\intR \big((\phi_{11}+\phi_{22})u_1u^{\star}_1+\phi_{12}u^{\star}_1 u_2\big)}{\intR u_1u^{\star}_1}.
\label{n4}
\eeq
Here we have introduced the notation $\phi_{ij}$ to denote the solution to the equation
 \[\Delta \phi_{ij}=\nabla\cdot\big(\frac{i}{2}[u_i\nabla u_j^*-u_j^*\nabla u_i]-u_i^*u_jhA_0\big)=0\quad\mbox{for}\;i,j=1,2\] subject
to homogeneous Neumann boundary conditions on $\partial\mathcal{R}$ and zero mean on $\mathcal{R}$,
so that $\tilde{\phi}(a^{\e}u_1+(a^{\e})^*u_2)$ in the nonlocal contribution $\mathcal{N}_2$ to the
nonlinearity $\mathcal{N}$ takes the form $\abs{a^{\e}}^2(\phi_{11}+\phi_{22})+(a^{\e})^2\phi_{21}+((a^{\e})^*)^2\phi_{12}.$

\begin{figure}
\begin{tabular}{c}
\begin{tabular}{cc}
\includegraphics[width=0.7\textwidth]{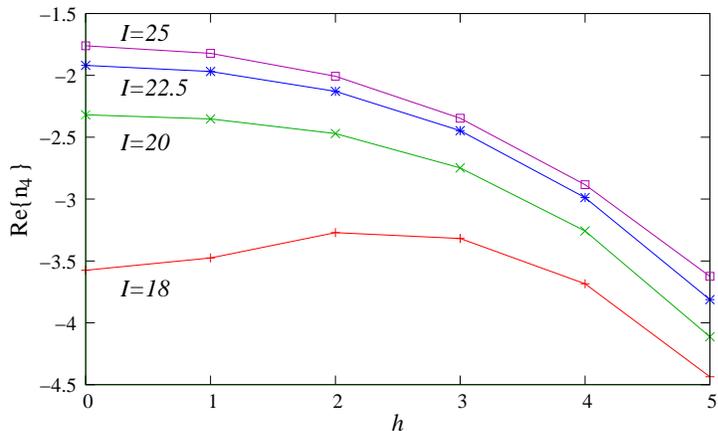}
\end{tabular}
\end{tabular}
\caption{Graph showing the required negativity of the quantity ${\rm{Re}}\,n_4$ as a function of $h$ for various values of $I$. Here we have taken $L=1$, $K=2/3$ and $\delta=4/15$. Other geometries were verified as well.} \label{n4graph}
\end{figure}

Provided ${\rm{Re}}\,n_4<0$, it is then easy to check from \eqref{aeqn} that the system undergoes a supercritical Hopf bifurcation to a periodic solution given by
\begin{equation}
a^{\e}(t)=\frac{\e^{1/2}}{\abs{{\rm{Re}}\,n_4}}e^{-i\big(\ipl+\gamma\e\big)t}\quad\mbox{where}\;\gamma:=\frac{{\rm{Im}}\,n_4}{{\rm{Re}}\,n_4}.\label{aepsform}
\end{equation}
We have verified the condition ${\rm{Re}}\,n_4<0$ numerically for a wide range of parameter values. See Figure \ref{n4graph}.

Summarizing the analysis above, we have shown:\\
\bthm \label{periodic} Fix a choice of parameters $K,L,\delta,h$ and $I$ such that the first eigenvalue of
$\mathcal{L}$ in \eqref{evpm} satisfies $\ipl\not=0$ and $n_4$ given by \eqref{n4} satisfies ${\rm{Re}}\,n_4<0$.
Then taking $\Gamma=\rpl+\e$ in
\eqref{psieqn}, there exists a value $\e_0>0$ such that for all positive $\e<\e_0$,
the system \eqref{psieqn}-\eqref{bc6} undergoes a supercritical Hopf bifurcation
to a periodic state $(\psi_{\e},\phi_{\e})$. Applying \eqref{verytang} to this solution, we see
that
\begin{equation}
\norm{\psi_{\e}-\bigg(a^{\e}(t)u_1+a^{\e}(t)^*u_1^{\dag}\bigg)}_{H^2(\mathcal{R})}\leq C\e^{3/2}\label{psiest}
\end{equation}
with $a^{\e}$ given by \eqref{aepsform}.
\ethm
\brk \label{stability} Though we do not present the analysis here, one can show that in fact this periodic solution
is asymptotically attracting. Since there is no unstable subspace associated with $\mathcal{L}_1$, nearby points
off of the center manifold $\mathcal{M}_{\e}$ are exponentially attracted to $\mathcal{M}_{\e}$. Then one argues
that nearby non-PT symmetric points on $\mathcal{M}_{\e}$ are attracted exponentially to the PT
symmetric part of $\mathcal{M}_{\e}.$ This type of argument was carried out in \cite{rsz} and no doubt the same type of argument will
work here though we have not checked the details.
\erk
\brk\label{stationary} One can also carry out the bifurcation analysis in the regime where $\lambda_1$ is real, which in particular
would correspond to parameter regimes where $I$ is sufficiently small.
As our primary goal in this article is to address issues related to periodic phenomena raised in \cite{agkns}, \cite{bepv} and \cite{bmp} we did not pursue
it here. In this case the center subspace is simply spanned by $u_1$ and a stable stationary state emerges for $\e>0$ of the form
\begin{equation}
\psi_{\e}\sim C\e^{1/2}u_1\label{statu}
\end{equation}
for some computable constant $C$. In the case where $h$ is positive
and sufficiently large, while $\lambda_1$ is real, one expects this stationary state to have magnetic vortices. Such a result might
be compared to the single vortex stationary solution found in \cite{dwz} for a similar model.
\erk

\section{Vortex formation}

\begin{figure}
\subfigure[]{\label{bbb1}
\includegraphics[width=6cm]{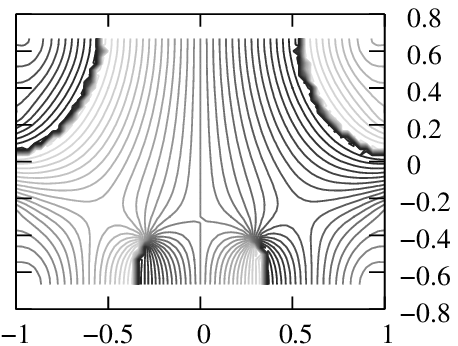}
}
\subfigure[]{\label{bbb2}
\includegraphics[width=6cm]{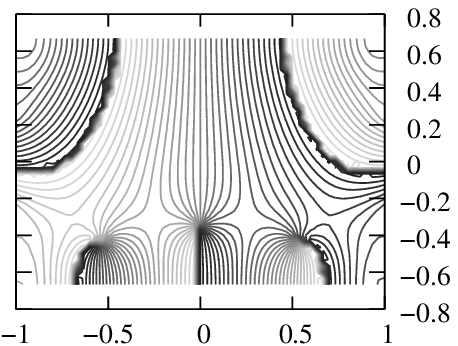}
}
\caption{Stationary magnetic vortices. The curves represent level sets of the phase of the first eigenfunction $u_1$. The geometry is $L=1, K=2/3, \delta = 4/15$. In (a) we used $h=15$, while in (b) we used $h=17$. In both cases we took $I=10$ which is below the critical value $I_c$, so that a stationary vortex solution emerges to the full problem
\eqref{psieqn}-\eqref{bc6} with leading profile given by $u_1$, cf. Remark \ref{stationary}}  \label{magnetic}
\end{figure}

As we mentioned in the Introduction, vortices form in this problem due to two separate effects. One type of vortex, that we term a magnetic vortex, is well-known.
 Magnetic vortices form as a result of the applied magnetic field, and can appear even when $I=0$. We present two examples of such vortices in Figure \ref{magnetic}. In both cases we used a rectangle  with parameters $L=1, K=2/3, \delta = 4/15$. In Figure \ref{magnetic}a the applied field is $h=15$, while in Figure \ref{magnetic}b the applied field is $h=17$. In both figures the current is $I=10$ which is below the critical current $I_c$ for this geometry and these values of $h$. Hence, the first eigenvalue of the linearized problem is real and so the center manifold here is one-dimensional. This is the regime discussed in Remark \ref{stationary}. There are other examples, not shown, where vortices form even when $I=0$.

What makes the present problem unusual is the formation of {\em kinematic} vortices, that is, vortices that are created even in the absence of magnetic fields. Physically, these vortices, just as the magnetic vortices, are points in space-time where the order parameter vanishes, the order parameter has a nonzero degree around these points, and large phase gradients occur near them. The formation of such kinematic vortices was extensively studied by Bendiyorov et al. \cite{bmp} using numerical simulations of the TDGL equations. The authors report on an unusual effect, where vortices appear periodically in pairs along the center line $x=0$ of the rectangle, and move along it. They also found that, depending on the parameter values in the problem, vortices can either form at opposite sides of the boundary and annihilate inside or else nucleate together at an internal point on the center line and move away from each other towards the boundary. When $h = 0$, both vortices appear a
 t the same time and they vanish, either by crossing the boundary simultaneously, or by annihilating each other, symmetrically about the line $y=0$. When $h \neq 0$ the $y$-symmetry is broken: the creation of these kinematic vortices can take place at different times, and their motion is not symmetric with respect to the line $y=0$.

We will use the theory developed in the preceding section to give a simple explanation for the formation and motion of kinematic vortices. Our analysis also gives a simpler means to compute when and where they form. In fact, we derive an explicit equation of motion for the kinematic vortices, fully based on the leading eigenfunction of the operator $ \mL$ defined in equation (\ref{evpm}). After deriving the equation of motion below, we demonstrate the different types of vortex creation and motion. One benefit of the new theory is that it allows us to easily detect additional types of kinematic vortex patterns, not observed in \cite{bmp}.

\begin{figure}
\subfigure[] {\label{ccc1}
\includegraphics[width=6cm]{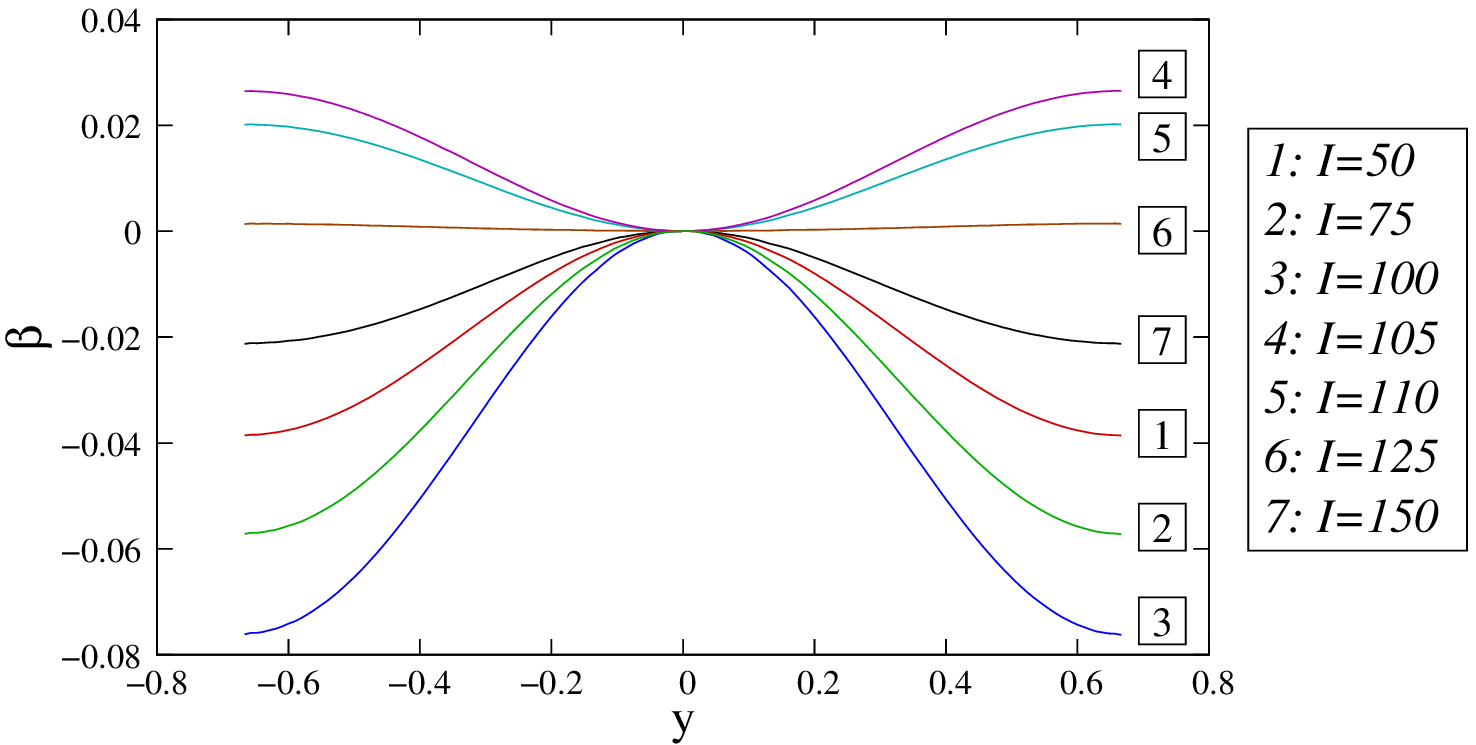}
}
\subfigure[] {\label{ccc2}
\includegraphics[width=6cm]{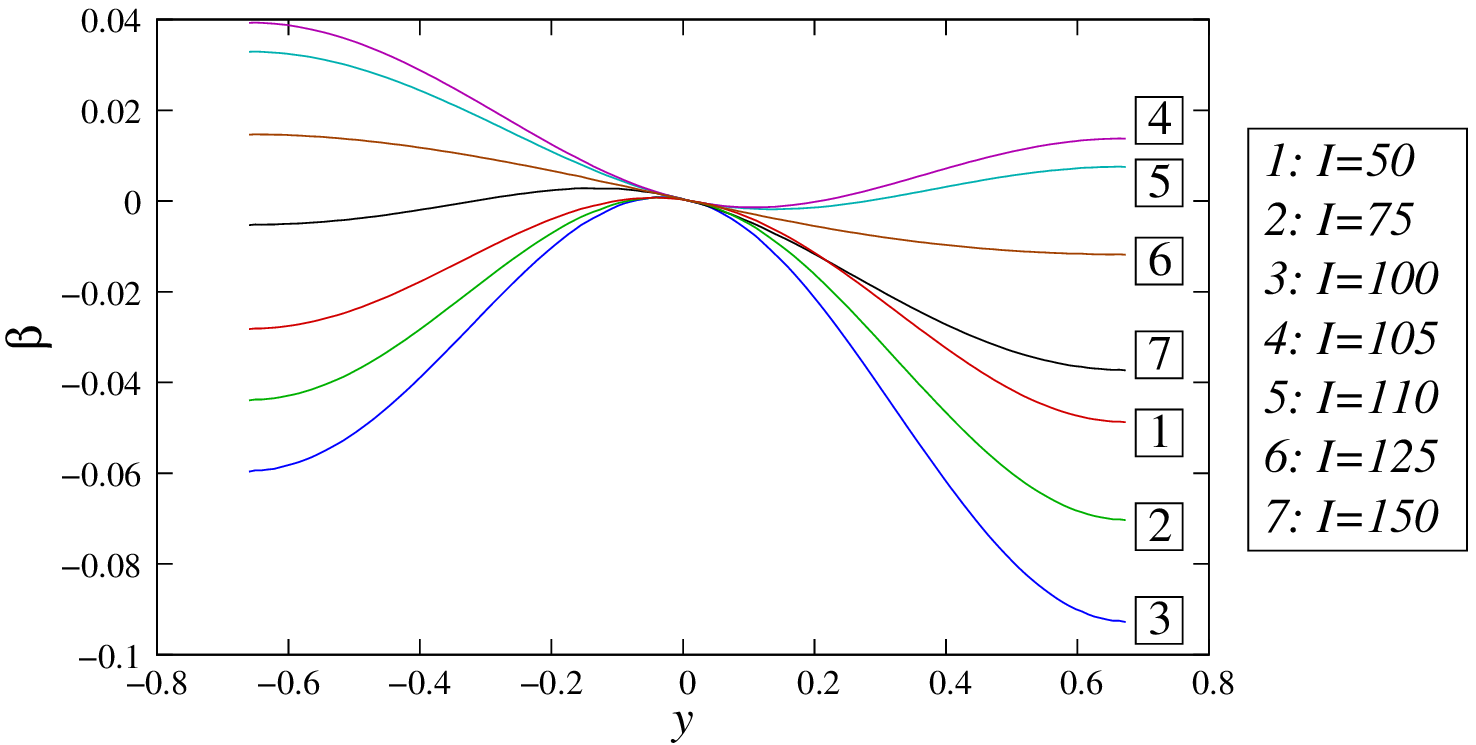}
}
\caption{Different shapes for the function $\beta(y)$. The geometry is $L=1, K=2/3, \delta = 4/15$. In (a) we set $h=0$, while in (b) we set $h=0.05$. } \label{varbeta}
\end{figure}

Mathematically, the formation and motion of the kinematic vortices are a consequence of the PT-symmetry of the problem. To see how they are created and move about, we consider the leading order term in the center manifold (cf. \eqref{psiest}): $$\psi = a^{\e}(t)u_1+a^{\e}(t)^*u_1^{\dag}.$$ It is convenient to introduce the notation $$a^{\e} = \xi \e^{-i \chi t},\;\; u_1(0,y) = g(y) e^{i \beta(y)},$$ where $\xi$ and $\chi$ take the values provided in equation (\ref{aepsform}). Therefore, along the rectangle's central line the order parameter is given (at leading order) by
\begin{equation}
\psi(0,y,t) = 2 \xi g(y) \cos\left(-\chi t + \beta(y)\right). \label{v1}
\end{equation}
Hence, the order parameter vanishes on the central line $x=0$ whenever the equation
\begin{equation}
\chi t = \beta(y) +\pi/2 + n \pi,\;\;  n=0, \pm 1, \pm 2, ...\label{v3}
\end{equation}
holds. Equation (\ref{v3}) is the equation of both motion {\em and} creation of kinematic vortices.
\brk\label{normal2} In the computation performed in this section we replaced the normalization condition (\ref{normal1}) with the normalization $u_1(0,0)=1$. This condition has the advantage that for all parameters we have $\beta(0)=0$, and therefore it is graphically easy to compare different  $\beta$ functions.
\erk
The boundary conditions on $u_1$ imply $\beta'(\pm K)=0$. Therefore, three simple shapes for $\beta$ might be expected: upward hump, downward hump, or a monotone shape. However, our numerical study shows that, while indeed  each of these shapes can occur, depending on the problem's parameters, other shapes are also present. Also, when $h=0$, the problem is symmetric with respect to the $y$ axis, and therefore $\beta$ is either an even function, as in Figure \ref{varbeta}a, or an odd function, if it is monotone. On the other hand, when the magnetic field is turned on, and $h \neq 0$, the symmetry of $\beta$ is broken.

To demonstrate some of the different possible shapes of $\beta(y)$ and their dependence on the parameters in the problem, we refer to Figure \ref{varbeta}. In all parts of this figure we took $L=1, K=2/3, \delta= 4/15$. In Figure \ref{varbeta}a we set $h=0$. Here $\beta$ appears as an even function of $y$.
Note that it changes its concavity depending on the values of the given current $I$.
In Figure \ref{varbeta}b we observe symmetry breaking in the shape of $\beta$ under the influence of a nonzero applied field $h=0.05$. In particular we point out that for some $I$ levels the function $\beta(y)$ becomes a monotone function. We are still pursuing an explanation for the dramatic change in $\beta$ as $I$ varies from 100 to 105.\\

We now analyze four examples to illustrate the effect of the shape of the function $\beta$ on vortex creation and motion.

\begin{figure}
\subfigure[]{\label{ddd1}
\includegraphics[width=6cm]{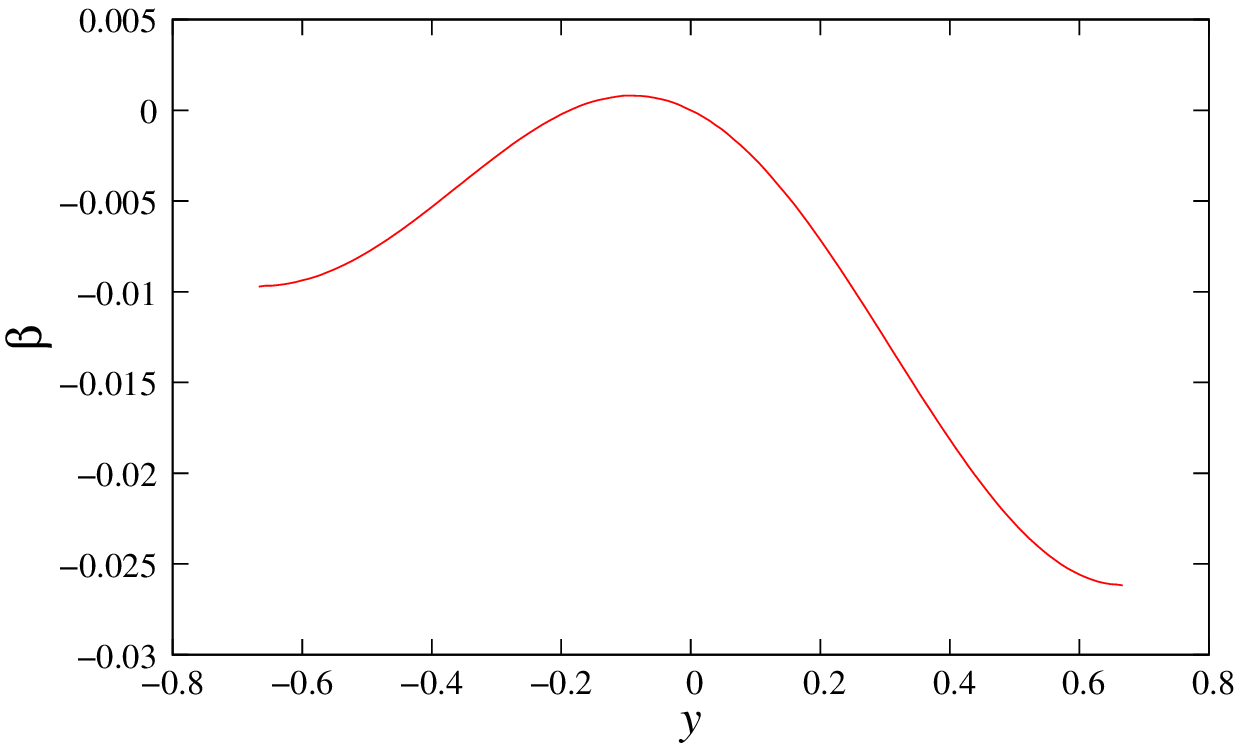}
}
\subfigure[]{\label{ddd2}
\includegraphics[width=6cm]{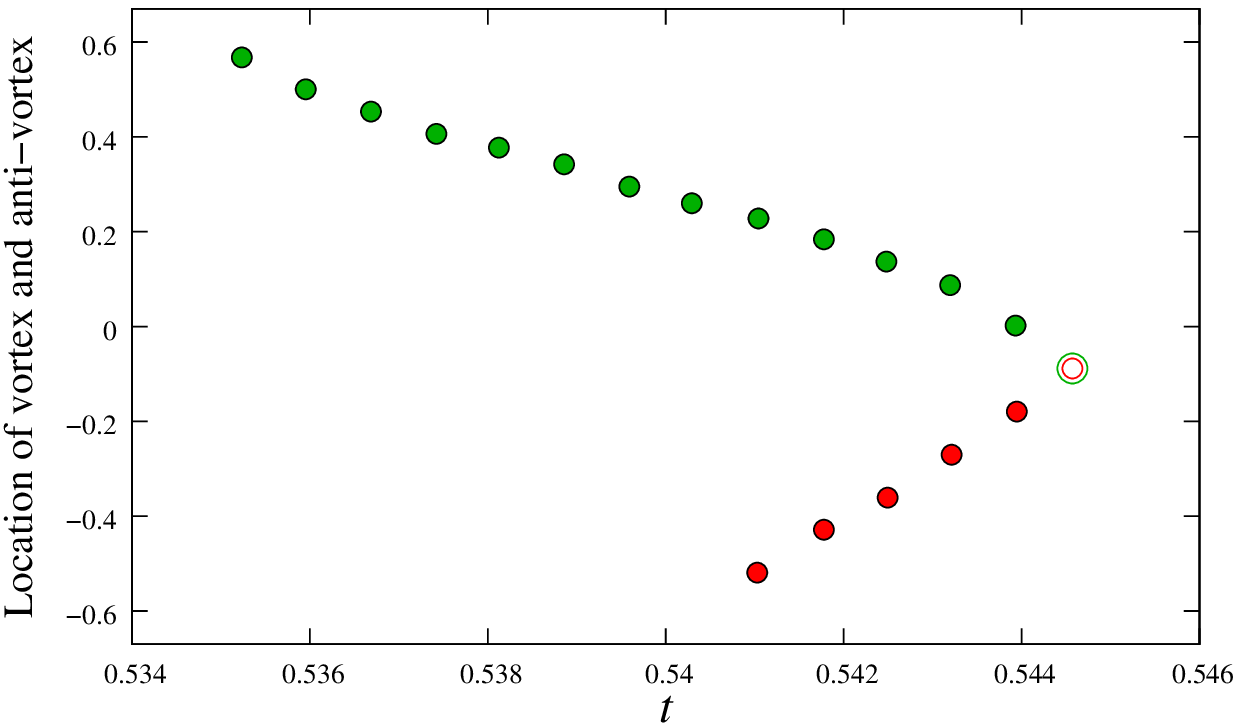}
}
\caption{Creation and motion of kinematic vortices. (a) The function $\beta(y)$ for the parameters $L=1, K=2/3, \delta = 4/15, h=0.05, I=25$. (b) The circles describe the location of the vortices in the $(y,t)$ plane.  } \label{beta1}
\end{figure}

\vskip 0.2cm \noindent {\bf Case 1:} We refer to Figure \ref{beta1}. We use the same geometric parameters as before, namely, $L=1, K=2/3, \delta= 4/15$. We also use $I=25, h=0.05$. The function $\beta$ is an asymmetric downward hump, as depicted in Figure \ref{beta1}a. The first solution of equation (\ref{v3}) occurs when $t$ is large enough so that $\chi t$ reaches $\beta(K) + \pi/2$
and a vortex emerges from the boundary. Then, as $t$ increases the solution moves to lower values of $y$. When $t$ is sufficiently large so that $\chi t = \beta(-K) + \pi/2$, a second vortex forms at the lower end $y=K$. As $t$ grows further, both vortices move towards each other. When $t$ reaches the level where $\chi t = \beta(y_m) + \pi/2$, where $y_m$ is the location of the maximal point of $\beta$, the vortices collide and they annihilate. This is an example of annihilation of a vortex/anti-vortex pair, where we use `anti-vortex' to refer to a vortex of negative degree.
Then, there is a time interval where equation (\ref{v3}) does not hold for any $y$, and therefore there are no kinematic vortices at those times. This scenario repeats itself when $\chi t = \beta (-K) + 3\pi/2$ and so on. The creation, motion and annihilation of vortices are shown in Figure \ref{beta1}b for a half of a single period.

\begin{figure}
\subfigure[]{\label{eee1}
\includegraphics[width=6cm]{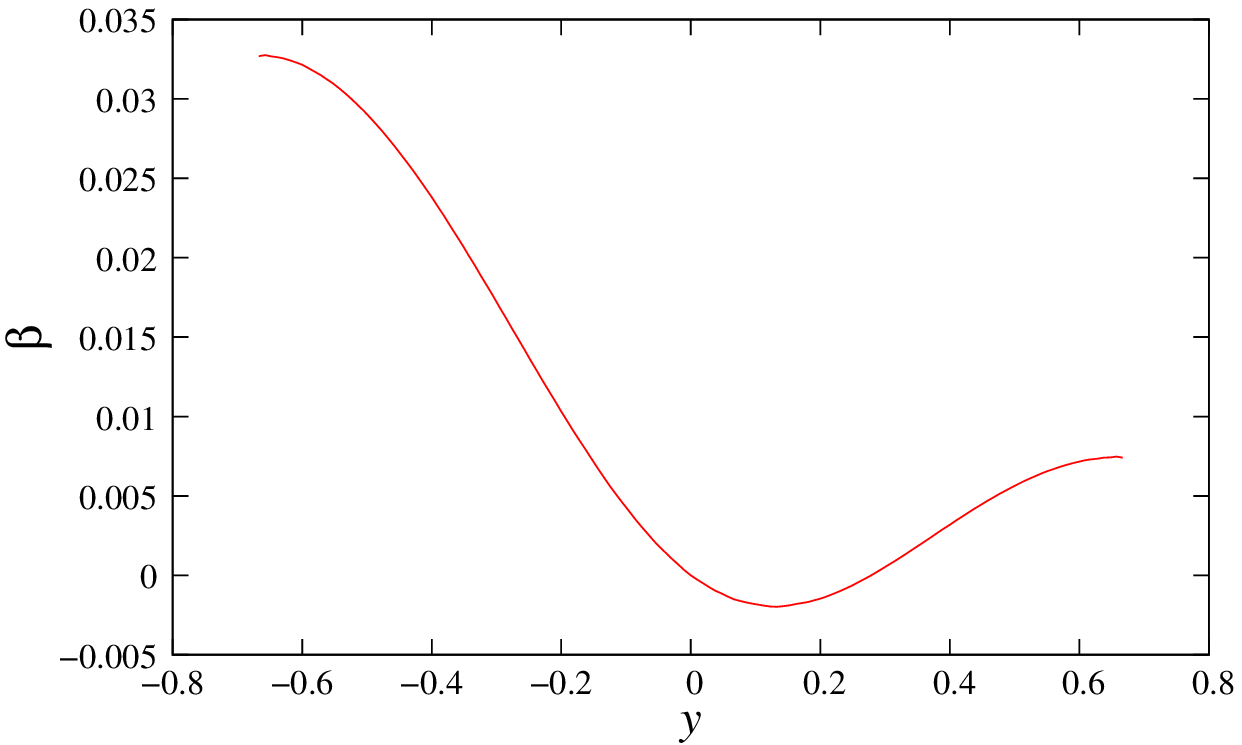}
}
\subfigure[]{\label{eee2}
\includegraphics[width=6cm]{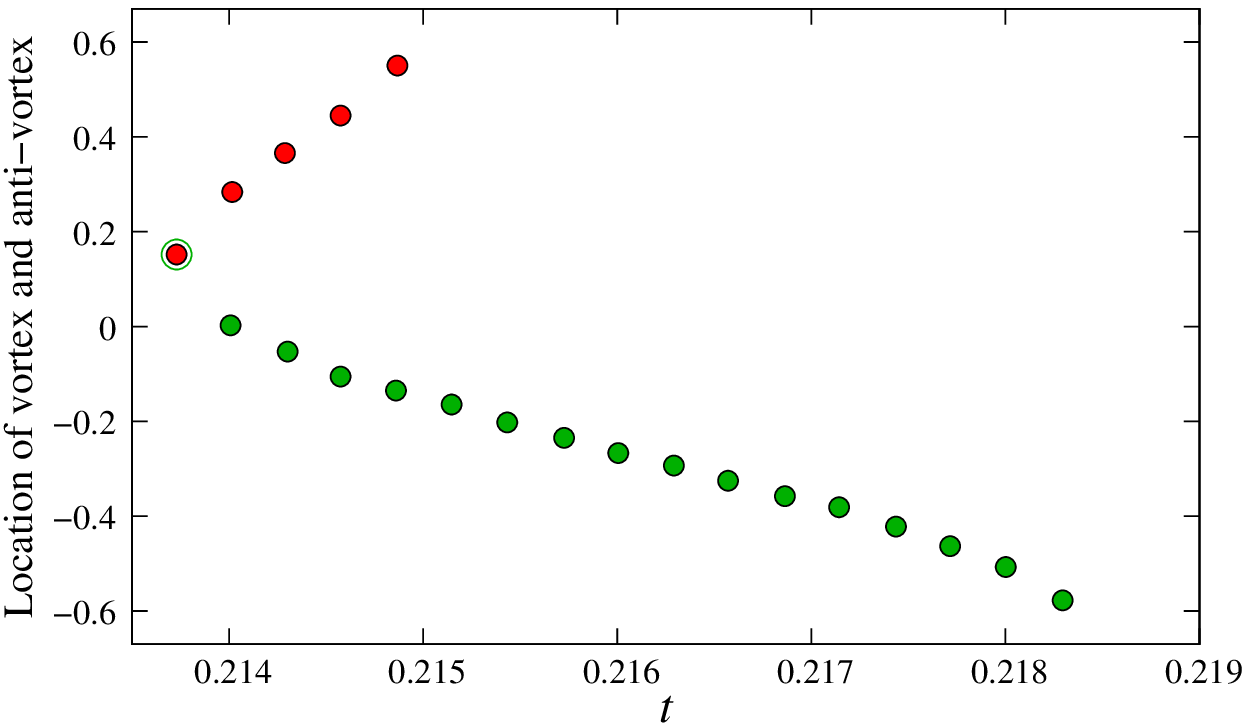}
}
\caption{ Creation and motion of kinematic vortices. (a) The function $\beta(y)$ for the parameters $L=1, K=2/3, \delta = 4/15, h=0.05, I=110$. (b) The circles describe the location of the vortices in the $(y,t)$ plane.  } \label{beta2}
\end{figure}

\vskip 0.1cm \noindent {\bf Case 2:} Using the same parameters as in the previous example, except increasing the current $I$ to take the value $I=110$, gives rise to a different shape for $\beta(y)$. As depicted in Figure \ref{beta2}a, it is now a distorted U shape. Therefore, denoting the location of the minimum of $\beta$ by $y_m$, when $t$ reaches the value where $\chi t = \beta(y_m) + \pi/2$, a vortex/anti-vortex pair is created inside the sample. As $t$ increases, equation (\ref{v3}) is satisfied at two locations, until a point of time where $\chi t = \beta(K) + \pi/2$. After that, only one vortex remains in the rectangle, and this vortex eventually leaves the domain when $\chi t = \beta(-K) + \pi/2$. The motion of these kinematic vortices is depicted in Figure \ref{beta2}b.

\begin{figure}
\subfigure[]{\label{ggg1}
\includegraphics[width=6cm]{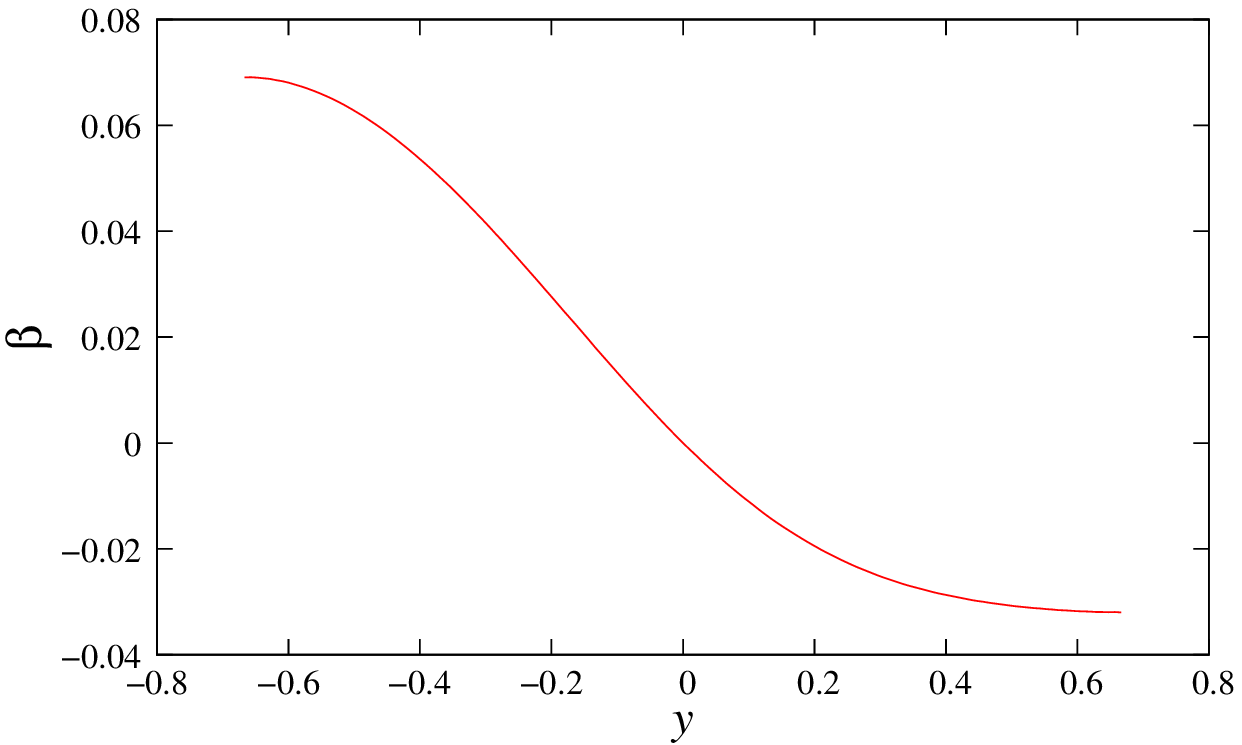}
}
\subfigure[]{\label{ggg2}
\includegraphics[width=6cm]{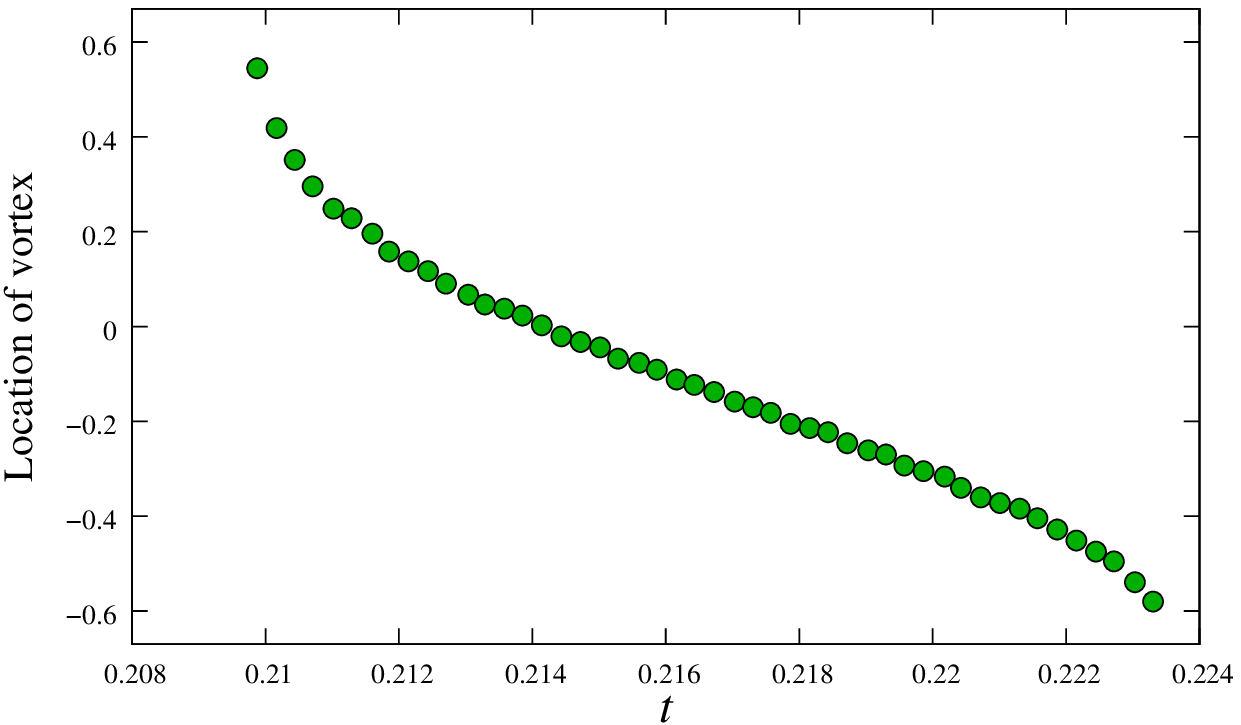}
}
\caption{ Creation and motion of kinematic vortices. (a) The function $\beta(y)$ for the parameters $L=1, K=2/3, \delta = 4/15, h=0.2, I=110$. (b) The circles describe the location of the vortices along the central line $x=0$ in the $(y,t)$ plane.  } \label{beta3}
\end{figure}

\vskip 0.1cm \noindent {\bf Case 3:} For the third example we maintain the same parameters as in the second example above, except that we increase $h$ to take the value $h=0.2$. For this choice of parameters $\beta(y)$ is a monotone function as depicted in Figure \ref{beta3}a. Now, the vortex is first created when $\chi t = \beta(K) + \pi/2$. It then travels to the lower end of the center line until $\chi t = \beta(-K) + \pi/2$. This motion is depicted in Figure \ref{beta3}b.

\begin{figure}
\subfigure[]{\label{kkk1}
\includegraphics[width=6cm]{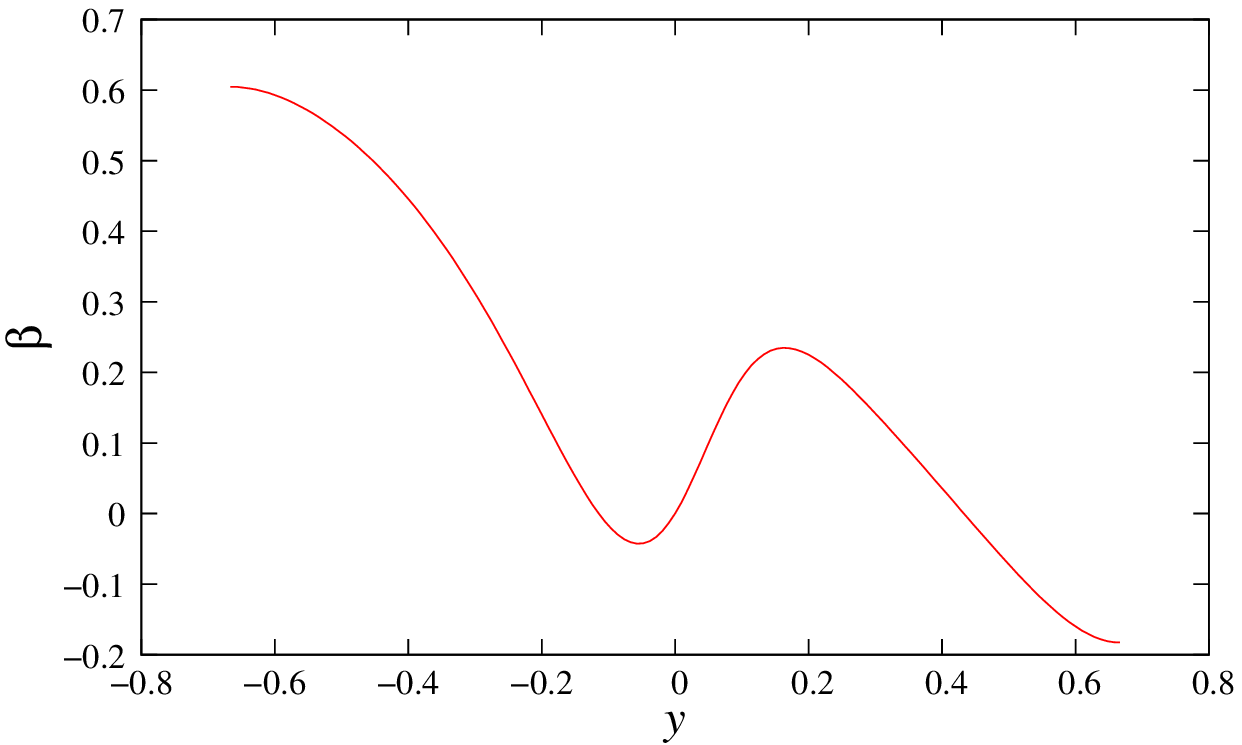}
}
\subfigure[]{\label{kkk2}
\includegraphics[width=6cm]{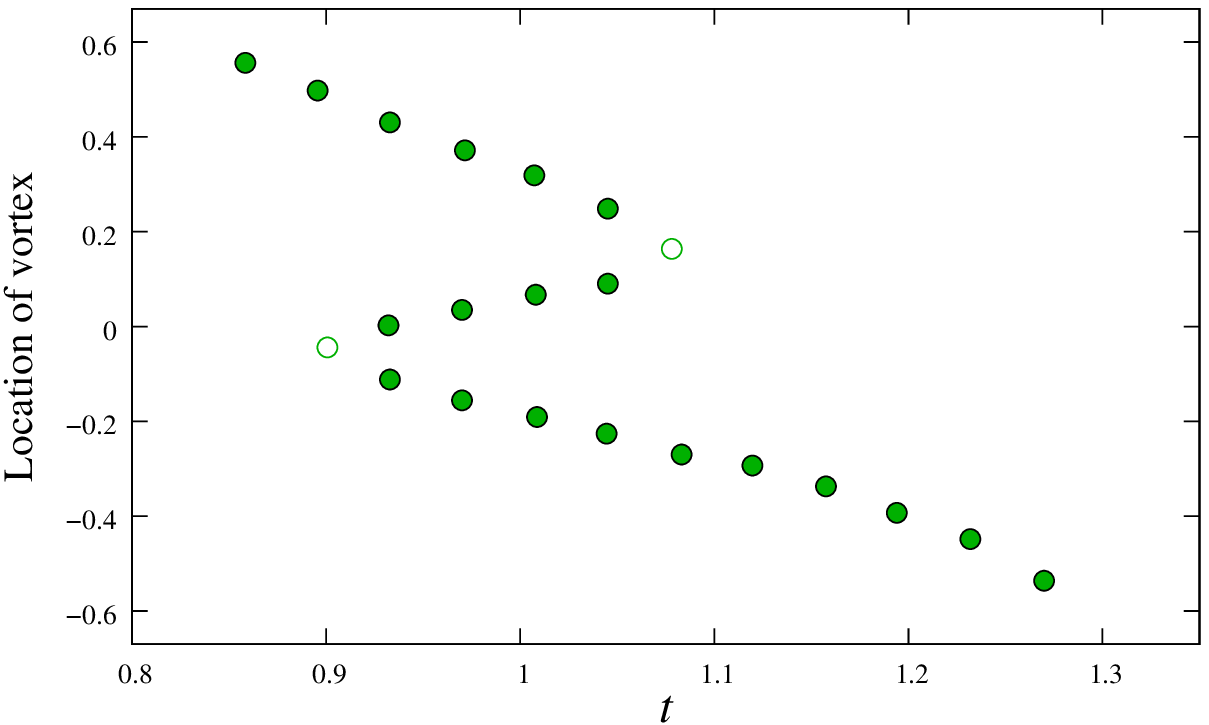}
}
\caption{ Creation and motion of kinematic vortices. (a) The function $\beta(y)$ for the parameters $L=1, K=2/3, \delta = 4/15, h=20, I=25$. (b) The circles describe the location of the vortices in the $(y,t)$ plane.  } \label{beta4}
\end{figure}

\vskip 0.1cm \noindent {\bf Case 4:} In the fourth example we present a case where $\beta(y)$ has both a local maximum {\em and} a local minimum. The function $\beta(y)$ for the parameters $h=20, \; I=25$ is shown in Figure \ref{beta4}a. This is the parameter choice of Figure \ref{var_beta1} as well. The geometry is the same as in the preceding examples. Following the scenarios above, if we look along the center line $x=0$ we expect to see a first vortex emerging at $y=K$ at a time that we denote $t_1$. Then, after this first vortex appears, a pair of vortices (of the same degree) appear at a later time, say $t_2$, where $\chi t_2=\beta(\bar{y})+\pi/2$ and $\bar{y}$ is the location of the interior local minimum of $\beta$. They move away from each other, until the one moving upward collides at a time $t_3$ with the first vortex moving downward.  Finally, the remaining vortex that moves downward reaches the boundary $y=-K$ and exits the rectangle. This vortex creation and moti
  on is indeed verified in Figure \ref{beta4}b.

\begin{figure}
\centering
\subfigure[]{\label{fff1} \includegraphics[width=6cm]{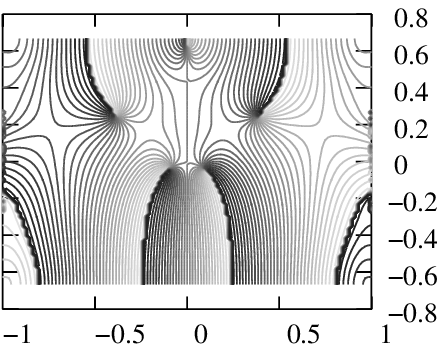}
}
\subfigure[]{\label{fff2}
\includegraphics[width=6cm]{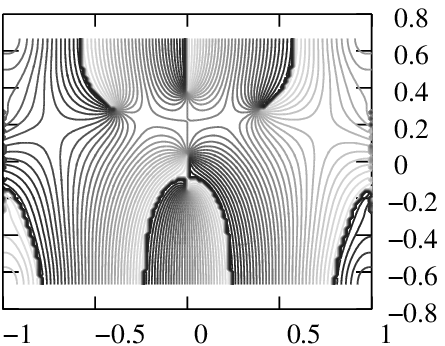}
}
\subfigure[]{\label{fff3}
\includegraphics[width=6cm]{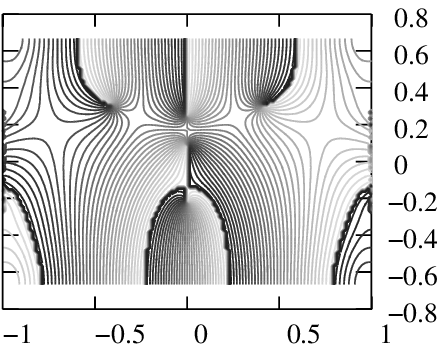}
}
\subfigure[]{\label{fff4}
\includegraphics[width=6cm]{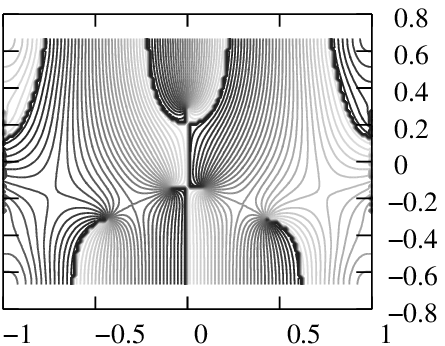}
}
\caption{Contour plots of the phase of the order parameter for the parameters of Case 4 at four distinct times during a period of the motion. In (a) we see five vortices, one of them at the upper end of the center line. In (b) two vortices that were earlier off the center line move into it and meet there. Then they separate and start moving away from each other along the center line in (c), with the upper one moving upward to approach a south-moving vortex along the center line. In (d) these two vortices seem to collide and then veer away from each other.
In these figures, the darkest lines are not significant in that they only represent a $-\pi$ to $\pi$ jump in the phase. The tips of these lines, however, represent the vortices. } \label{5vortex_1}
\end{figure}

This example, however, has several peculiar features.  In Cases 1 and 2 vortices formed or disappeared in pairs of vortex-antivortex structure on the center line. In Case 4, on the other hand, the picture is different. We refer to a sequence of snapshots in Figure \ref{5vortex_1}. In Figure \ref{5vortex_1}a we observe five vortices, with only one on the center line. The two vortices far from the center line barely move throughout the period of the evolution. We view these two ``sluggish" vortices as magnetic vortices. The vortex {\em on} the center line is the kinematic vortex that formed at $y=K$ at $t=t_1$ as explained above. The two vortices located on either side of the center line, and near it, which are of the same degree, move towards each other. Eventually they meet at $t=t_2$ on the center line, giving rise to the two kinematic vortices that were discussed above, and are shown in Figure \ref{beta4}b as well as Figure \ref{5vortex_1}b. These two vortices move as descr
 ibed above until $t=t_3$. Then, the middle vortex on the center line meets the upper vortex on the center line. This is shown also in Figure \ref{5vortex_1}c. Then, as shown in Figure \ref{5vortex_1}d, this pair of vortices split away from the center line. This new pair of vortices  moves away from the center line and upward, and the entire process repeats itself periodically.

The picture we just outlined indicates that kinematic vortices can move away from the center line. We therefore also term such vortices kinematic, namely those
that spend part of a period on and part of a period off the center line, since their formation and motion follow directly from the PT symmetry and the structure it imposes on the center manifold.

\section{Discussion}

All of the examples from the previous section illustrate that, at least near the normal state, there is a dichotomy in vortex behavior when both applied currents and applied magnetic fields are present. The expansion based on center manifold reduction gives a partial explanation for this phenomenon, with the kinematic vortices arising in part due to the PT symmetry of the problem. In any event, it is clear that the variety of possible vortex behavior in this system is far more extensive than that seen in models capturing only magnetic effects. In particular, the motion law \eqref{v3}, based on small amplitude asymptotics rather than large Ginzburg-Landau parameter asymptotics as is more common in the literature, allows for a wide range of effects including boundary and interior nucleation, collision of
like-signed vortices and periodicity of these events. Of course, all of the rigorous analysis we conduct necessarily involves small amplitude solutions since it is based on a bifurcation from the normal state. One would imagine that an even richer array of vortex behavior is possible for this system if one looks far from the normal state though an analytical approach would clearly require different tools.
\newpage
\noindent
{\bf Acknowledgments.} L. Peres Hari and J. Rubinstein were generously supported by an ISF grant. P. Sternberg was generously supported by NSF grant DMS-1101290 and a Simons Foundation Collaboration Grant.

\end{document}